\documentclass[10pt,a4paper, english]{article}
\usepackage{amssymb}
\usepackage{fontenc}
\usepackage{amsthm}
\usepackage{amsmath}
\usepackage[colorlinks=false]{hyperref}
\newtheorem{teo}{Theorem}[section]
\newtheorem{defi}[teo]{Definition}
\newtheorem*{eje}{Example}
\newtheorem*{nota}{Remark}
\newtheorem{lema}[teo]{Lemma}

\newtheorem{coro}[teo]{Corollary}
\newtheorem{prop}[teo]{Proposition}

\usepackage{graphicx}
\usepackage[dvips]{pstcol}
\usepackage{epsfig}
\usepackage{fancyhdr}

\newenvironment{demo}{
  {\noindent\sc{Proof:}}\
  }{$\fbox{}$ \par \bigskip}

\begin{document}
\title{On dilatation factors of braids on three strands}
\author{Marta Aguilera\footnote{Partially supported by MTM2010-19355, P09-FQM-5112 and FEDER.}\\ Universidad de Sevilla}
\date{}
\maketitle
\begin{abstract}
In this work we present a natural surjective map from rigid braids in $B_3$ (in Garside sense) to $SL_2(\mathbb{N})$. This map provides an upper and a lower bound for the dilatation factor of a pseudo-Anosov $3$-strand braid. These bounds only depend on the canonical length of the classical Garside structure of $B_3$.
\end{abstract}

\section{Introduction}

\begin{figure}[h!]
\centerline{\includegraphics*{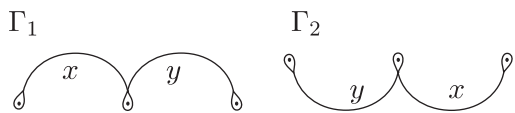}}
\caption{Train tracks for rigid braids on three strands.}
\label{2tt}
\end{figure}

In this paper we review some well-known results about the braid group in three strands, and we rewrite them in terms of the Garside structure. In this way, we see that the dynamic of a braid in a super summit set (i.e. with minimal length in its conjugacy class) is easy to describe.
If such a braid is reducible, it must be $\Delta^{2s} \sigma_1^k$ or $\Delta^{2s}\sigma_2^k$, $s,k \in \mathbb{Z}$. Their reduction systems are simple: a curve around the first two punctures for the former, and a curve around the last two punctures for the latter. 
 In the pseudo-Anosov case, we will point out that, despite the fact that there are only two train track graphs $\Gamma_1, \Gamma_2$ as in Figure \ref{2tt} that carry every foliation in the 3-times punctured disc $\mathbb{D}_3$, these are not train tracks for every pseudo-Anosov braid. However, if a braid $\beta \in B_3$ is in a super summit set, it is rigid (a Garside theoretical property), and either $\Gamma_1$ or $\Gamma_2$ is a train track for $\beta$ (\cite{Han}, \cite{Mura}).

In addition, if $\beta$ is in a super summit set, it is straightforward to obtain the associated matrix, and therefore the foliation and dilatation factor, from its Garside normal form. This allows us to study the dynamic of any braid through a conjugate in the super summit set. This also gives a natural map from rigid braids to $2\times2$ matrices with non-negative entries:
\begin{teo}\label{corointro}
There exists a $2$-to-$1$ surjective map from the set of rigid $3$-braids (modulo $\Delta^2$) to $SL_2(\mathbb{N})$, which sends each rigid braid to the matrix associated to its corresponding train track $\Gamma_i$, $i=1,2$.
\end{teo}

As a consequence, given a pseudo-Anosov braid with minimal length $\ell$ in its conjugacy class, we can give a lower and an upper bound for its dilatation factor:

\begin{teo}\label{teo1}
Let $\beta \in B_3$ be a pseudo-Anosov braid in its super summit set. If $\beta$ has canonical length $\ell$ and dilatation factor $\lambda$, then:
\[\frac{1}{2} \left( \ell +1 + \sqrt{ \left( \ell+3 \right)  \left( \ell-1  \right) } \right)\leq \lambda,\]
and:
\[\begin{array}{lcl}
    \lambda \leq \phi^\ell & \quad & \mbox{ if } \ell \mbox{ is even,} \\
    \lambda \leq F_{\ell}+\sqrt{F_{\ell}^2-1}< \frac{2}{\sqrt{5}}\left( \phi^\ell+\phi^{-\ell}\right) & \quad & \mbox{ if } \ell \mbox{ is odd,} \\
\end{array}\]
where $\phi$ is the golden ratio, and $F_i$ is the $i^{th}$ Fibonacci number.

The lower bound is a minimum, reached only by the conjugates of $\beta=\sigma_1\sigma_2^{-(\ell-1)}$.
The upper bound is a maximum if $\ell$ is even, reached only by the conjugates of $\beta=\left(\sigma_1\sigma_2^{-1}\right)^{\frac{\ell}{2}}$. If $\ell$ is odd, the maximal dilatation factor reached only by the conjugates of $\beta=\left(\sigma_1\sigma_2^{-1}\right)^{\frac{\ell-1}{2}}\sigma_1$ satisfies that $2F_\ell-1 < \lambda < 2F_\ell$ with $2F_\ell=\frac{2}{\sqrt{5}}\left( \phi^\ell+\phi^{-\ell}\right)$.
\end{teo}

\noindent \textbf{Acknowledgements.} I would like to thank Dan Margalit, Bert Wiest, J\'er\^{o}me Los, Eiko Kin and Spencer Dowdall for helpful discussions and comments, some of them during a stay at Centre de Recerca Mat\`ematica, to whom I also want to thank for its hospitality. I am also grateful to my advisor Juan Gonz\'alez-Meneses for numerous useful comments, corrections and suggestions on preliminary drafts of this paper.

\section{Background}
The braid group on $n$ strands $B_n$ is isomorphic to the mapping class group of the $n$ times punctured disc fixing the boundary pointwise $MCG(\mathbb{D}_n, \partial\mathbb{D}_n)$. Collapsing the boundary of $\mathbb{D}_n$ to a point, the mapping class group of the $n$-times punctured disc can be considered as a subgroup of the $(n+1)$ times punctured sphere: it will be denoted $MCG(\mathbb{D}_n)$.

This geometric approach to $B_n$ allows the use of the Nielsen-Thurston classification theorem \cite{FLP}, \cite{Thurston}, \cite{Thur.notes}. A braid $\beta$ is periodic if there exist $k\in \mathbb{N}, s\in \mathbb{Z}$, such that $\beta^k=\Delta^{2s}$, where $\Delta^2$ is a Dehn twist along the disc's boundary $\partial\mathbb{D}_n$. $\beta$ is reducible if there exists a non-degenerate 1-manifold $C\subset\mathbb{D}_n$ fixed by $\beta$. Finally, $\beta$ is pseudo-Anosov if there exists a pair of transverse measured foliations  $(\mathcal{F}_s, \mu_s)$ y $(\mathcal{F}_u,\mu_u)$, and a real number $\lambda >1$, such that $\beta \left( (\mathcal{F}_s, \mu_s)\right)=(\mathcal{F}_s, \lambda\mu_s)$ and $\beta((\mathcal{F}_u, \mu_u))=(\mathcal{F}_u, \lambda^{-1}\mu_u)$.  The classes of reducible and periodic braids are not disjoint, so from now on we will call reducible those elements which are reducible and non-periodic.

The classification problem can be solved using the train tracks techniques introduced by Bestvina-Handel in the nineties \cite{BH}, \cite{BH2}. These are combinatorial objects which encode the dynamics on the surface in terms of linear algebra. In the pseudo-Anosov case, they also give the structure of the unstable foliation $\mathcal{F}^u$, its measure and the dilatation factor $\lambda$.

The dynamic of a braid only depends on its conjugacy class, in particular so does its Nielsen-Thurston type. So we can study the geometry of any given $\beta\in B_n$, through any conjugate $\widetilde{\beta}=\alpha^{-1}\beta\alpha$: Periodicity is easily recognizable in braid groups \cite{BGGM3}, and $\beta^k=\Delta^{2s}$ if and only if $\widetilde{\beta}^k=\Delta^{2s}$. If $\beta$ is reducible, and $C$ is a \textit{reduction system} for $\beta$ (1-manifold such that $\beta(C)=C$), then $\alpha(C)$ is a reduction system for $\widetilde{\beta}$. In the pseudo-Anosov case, if $\Gamma$ is a train track for $\beta$, then $\alpha(\Gamma)$ is train track for $\widetilde{\beta}$. Also the combinatorial maps associated to $\Gamma$ and $\alpha(\Gamma)$ are the same $\widetilde{f}_{\beta,\Gamma}=\widetilde{f}_{\widetilde{\beta},\alpha(\Gamma)}$, and so are the matrices $M(\beta,\Gamma)=M(\widetilde{\beta}, \alpha(\Gamma))$ (see next section for definitions).

\subsection{Train tracks}

In this section we will review some basic facts about train tracks (\cite{BH}, \cite{BH2}, \cite{F-M}).

Let $\Gamma=(V,E)$ be a labeled graph embedded into the punctured disc $\mathbb{D}_n$, such that each component of $\mathbb{D}_n \setminus \Gamma$ is a punctured disk or a ring and $\pi_1(\Gamma)=\pi_1(\mathbb{D}_n)=F_n$.
 In the class of a given $\beta\in MCG(\mathbb{D}_n)$ there exists a representative automorphism $f_{\beta}$ which maps $\Gamma$ into a tubular neighborhood $U$ of itself. The composition of this map with a deformation retract $\iota: U \rightarrow \Gamma$, allows us to associate to each $\beta$, a map $\widetilde{f}_\beta:\Gamma\rightarrow \Gamma$. Notice that we can chose $f_\beta$ and $\iota$ such that  $\widetilde{f_\beta}(V)\subseteq V$, and such that $\iota \circ f_\beta$ is injective in the interior of each edge, so for any $e\in E$, $\widetilde{f}_\beta(e)$ is an edge path. Thus $\widetilde{f}_\beta$ can be seen as a combinatorial map.

Removing vertices of valence $2$, and contracting edges which end in a valence $1$ vertex, we can suppose such a graph $\Gamma$ to have all vertices of valence at least three. In addition, we will assume that at each vertex there is a well defined tangent, so we can distinguish between those edges entering from one direction and those entering from the other. The labels of the edges must satisfy the switch condition: the sum of the labels going in from one side must be equal to the sum of those going out.

A combinatorial map $\widetilde{f}$ \textit{backtracks} if there is an edge $a\in E$, and $k>0$ such that $\widetilde{f}^k(a)$ contains the subword $ee^{-1}$ or $e^{-1}e$, for some $e\in E$.
A combinatorial map is said to be \textit{efficient} if it does not backtrack. Notice that $\widetilde{f}_\beta$ is efficient as combinatorial map if $\left(\iota \circ f_\beta\right)^k$ is injective in the interior of all edges for all $k>0$.

Given $\beta\in MCG(\mathbb{D}_n)$, the graph $\Gamma$, with the properties described above, is a \textit{train track} graph for $\beta$, if the combinatorial map $\widetilde{f}_\beta^k$ is efficient $\forall k\geq1$. In this case, the tangencies at the vertices of $\Gamma$ can be chosen so that $\iota \circ f_{\beta}$ respect tangencies, that is the image of each edge is a smooth edge path.

We can associate to a braid $\beta$ and $\Gamma$ a \textit{transition matrix} $M=M(\beta,\Gamma)$, where each entry $m_{i,j}$ is the number of times the edge $e_j$ appears in $\widetilde{f_\beta}(e_i)$. Obviously, $M$ has non-negative entries. If the graph $\Gamma$ is a train track for $\beta$, $M(\beta, \Gamma)$ contains the geometric information about $\beta$.

If the transition matrix $M$ is reducible (i.e. $M^k$ has at least one zero entry for all $k\in \mathbb{N}$), then the element $\beta$ is reducible. If the matrix is not reducible, then $\beta$ is either periodic or pseudo-Anosov. The Perron-Frobenius theorem states that, in the irreducible case, the greatest eigenvalue $\lambda$ is real, has multiplicity $1$ and $\lambda\geq 1$. The element $\beta$ is periodic if and only if $\lambda=1$, and it is pseudo-Anosov if and only if $\lambda>1$. In the latter case, the eigenvalue $\lambda$ coincides with the dilatation factor, and the associated eigenvector yields a measure for $\Gamma$ that encodes the unstable foliation $\mathcal{F}^u$ \cite{FLP}.

\subsection{Garside structure}
Braid groups have a well-known presentation \cite{Bir}:
\[B_n=\left\langle \sigma_1,\sigma_2,\ldots, \sigma_{n-1} : \quad
\begin{array}{cl}
  \sigma_i\sigma_j=\sigma_j\sigma_i & |i-j|>1 \\
  \sigma_i\sigma_j\sigma_i=\sigma_j\sigma_i\sigma_j & |i-j|=1
\end{array}
\right\rangle .
\]

\begin{figure}[h!]
\centerline{\includegraphics*{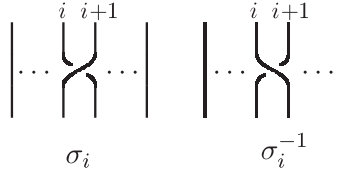}}
\caption{Positive and negative crossings in $B_n$.}
\label{crossings}
\end{figure}

These groups can be endowed with the classic Garside structure \cite{DehPar}, \cite{Gar}, that is a triple $(B_n, B_n^+,\Delta)$, where $B_n^+$ is the monoid generated by the \textit{positive crossings} $\sigma_1,\ldots,\sigma_{n-1}$ (see Figure \ref{crossings}), and $\Delta\in B_n^+$ is called the Garside element, \[\Delta=\left(\sigma_1\cdots\sigma_{n-1}\right)\left(\sigma_1\cdots\sigma_{n-2}\right)\cdots\left(\sigma_1\sigma_2\right)\sigma_1.\]
 For any $\gamma, \ \beta\in B_n$, we will say that $\gamma$ is a prefix of $\beta$, $\gamma\preccurlyeq\beta$, if $\gamma^{-1}\beta\in B_n^+$. This is a partial order that endows $B_n$ with a lattice structure (with well-defined $gcd$ and $lcm$), used to define normal forms.

The mapping class $\Delta$ is represented by a half Dehn twist along the disc's boundary. Considering braids up to Dehn twists along the boundary is equivalent to collapse the boundary to one point, so $MCG(\mathbb{D}_n)=B_n/\langle\Delta^2\rangle$. From the geometric point of view $\beta$ has the same properties as $\Delta^{2i}\beta$, $\forall i\in \mathbb{Z}$, so we will usually consider braids up to multiplication by $\Delta^2$. It is easy to see that $\Delta^2$ is central, actually it generates the center  $Z(B_n)=\langle\Delta^2\rangle$, \cite{chow}.

The Garside element $\Delta$ also satisfies:
\begin{itemize}
  \item $[1,\Delta]=\{s\in B_n \, :\, 1 \preccurlyeq s \preccurlyeq \Delta\}$, the set of \textit{simple elements}, generates $B_n$.
  \item The inner automorphism $\tau$ corresponding to conjugation by $\Delta$ preserves the lattice structure. Equivalently:
      \[\tau(B_n^+)=\Delta^{-1}B_n^+ \Delta=B_n^+.\]
\end{itemize}

\begin{defi}{\rm \cite{EM}}
Given a braid $\beta\in B_n$, the decomposition $\beta=\Delta^s s_1 s_2 \cdots s_k$ is its left-normal form if $1\prec s_i \prec \Delta \ \forall i$ and $s_i$ is the greatest simple prefix of $s_i s_{i+1}\cdots s_k$, $\forall i=1\ldots k$. The integer $k$ is the canonical length of $\beta$.
\end{defi}

The computational complexity of the calculation of the normal form of a braid on $n$ strands written as a product of $\ell$ generators $\sigma_i^{\pm1}$ is $O(\ell^2n\log n)$ for $n>3$ \cite{WP}, and linear in $\ell$ if $n=3$.

\begin{defi}{\rm \cite{BGGM}}
Let $\beta$ be a braid and $\beta=\Delta^ss_1\cdots s_k$ its normal form. We will say the braid $\beta$ is rigid if $k\geq1$ and the product $s_k\tau^{-s}(s_1)$ is in left normal form, or if $\beta=\Delta^s$, $s\in \mathbb{Z}$.
\end{defi}

If a braid $\beta$ has normal form as above, then $\tau(\beta)=\Delta^s\tau(s_1)\ldots \tau(s_k)$ is in normal form. Therefore, $\beta$ is rigid if and only if $\tau(\beta)$ is rigid.

\begin{teo}{\rm \cite{BGGM}}
For every pseudo-Anosov $\beta\in B_n$, there exists $0<k\leq\binom{n}{2}$ such that $\beta^{k}$ is conjugate to a rigid braid. In $B_3$, one can take $k=1$, that is, every pseudo-Anosov braid is conjugate to a rigid braid.
\end{teo}

In the braid group $B_3$, we will see how to extract geometric information of a pseudo-Anosov braid via a rigid conjugate.

\section{Braid group in three strands}
\subsection{Normal forms in $B_3$}
The braid group in three strands has specially nice properties. Normal forms are easily computable and rigidity is directly recognizable.
The Garside element is $\Delta=\sigma_1\sigma_2\sigma_1=\sigma_2\sigma_1\sigma_2$, and the simple elements are $\{1, \sigma_1, \sigma_2, \sigma_1\sigma_2, \sigma_2\sigma_1, \Delta\}$.
Given a braid as a concatenation of $\sigma_i^{\pm1}$, it is easy to rewrite it as the product of $\Delta^s$, $s\in \mathbb{Z}$, and a positive word in $\{\sigma_1, \sigma_2\}$, using the equalities:
 \[\begin{array}{ccc}
    \sigma_1^{-1}=\Delta^{-1}\sigma_1\sigma_2&\quad & \sigma_1\Delta^{\pm1}=\Delta^{\pm1}\sigma_2\\
    \sigma_2^{-1}=\Delta^{-1}\sigma_2\sigma_1&\quad & \sigma_2\Delta^{\pm1}=\Delta^{\pm1}\sigma_1\\ \end{array}\]

Notice that each simple element, except $\Delta$, can be written in a unique way as a word in $\{\sigma_1, \sigma_2\}$. The product of two proper simple elements $s_1s_2$ is in normal form if and only if the last letter of $s_1$ equals the first letter of $s_2$. Therefore, if $s_is_{i+1}$ is in normal form, the last crossing in $s_i$ and the length of $s_{i+1}$ characterize the factor $s_{i+1}$. Repeating this process, a product of proper simple elements in normal form $s_1\cdots s_k$ is determined by giving $s_1$ and the length of the factors $s_2, \ldots, s_k$.

So, the normal form of any braid $\beta=\Delta^s s_1\cdots s_k$, $k>0$, can be codified by the tuple $(s; i; p_1, q_1,\ldots, p_r, q_r)\in \mathbb{Z}\times \{1,2\} \times \mathbb{N}^{2r}$, $r>0$. The integer $s$ is the exponent of $\Delta$. The first crossing in $s_1$ is $\sigma_i$. The other elements in the tuple indicate the length of the factors $s_1,\ldots, s_k$ in the following way. The first $p_1$ simple elements in the normal form have length one (each one consists of one crossing). Then, they are followed by $q_1$ elements of length two, and after those there are $p_2$ elements of length one, etc. That is, for any $j=1, \ldots, r$, the simple factors from position $\sum_{i<j}\left(p_i+q_i\right)+1$ to position $\sum_{i<j}\left(p_i+q_i\right)+p_j$ have length one, and those from position $\sum_{i<j}\left(p_i+q_i\right)+p_j+1$ to position $\sum_{i<j+1}\left(p_i+q_i\right)$ have length two.
For a coherent notation only $p_1$ and $q_r$ could be zero. We will codify the braid $\Delta^s$ by $(s;1;0,0)=(s;2;0,0)$.

\begin{eje}
$(-4; 1; 0, 1,3,2)$ corresponds to the braid
\[\beta=\Delta^{-4}\sigma_1\sigma_2 \cdot \sigma_2 \cdot \sigma_2 \cdot \sigma_2 \cdot \sigma_2\sigma_1 \cdot \sigma_1\sigma_2.\]
\end{eje}

Recall that a braid $\beta\neq\Delta^s$, with normal form $\beta=\Delta^ss_1\cdots s_k$ is rigid if $s_k\tau^{-s}(s_1)$ is in normal form. That is, the normal form of $\Delta^{-s}\beta$ must start and finish with the same letter $\sigma_i$ if $s$ is even, and it must start and finish with different letters if $s$ is odd. The reader can check that a braid $\beta$ is rigid, in terms of the associated tuple, if and only if $s +\sum q_j$ is even. It is easy to check that $\tau\left((s;1;p_1,q_1,\ldots,p_r,q_r)\right)=(s;2;p_1,q_1,\ldots,p_r,q_r)$. If a braid $\beta=\Delta^s s_1\cdots s_k$ with $k>1$ is not rigid, then its conjugate by $\tau^{-s}(s_1)$ either is rigid, or its canonical length is strictly smaller than $k$ (or both things happen). This conjugation is known as \textit{cycling} \cite{EM}, and a finite number of iterations provides a rigid conjugate of any initial non-periodic braid. The Nielsen-Thurston type of a rigid braid it is easily recognizable:

\begin{prop}{\rm \cite{Mura}}
Every braid $\beta\in B_3$ is conjugate to a braid $\alpha=(s;1;p_1,q_1,\ldots, p_r, q_r)$ called Murasugi representative of $\beta$ such that:
\begin{itemize}
  \item If $\beta$ is periodic, $r=1$ and \begin{itemize}
                                            \item $s\in \mathbb{Z}$ and $(p_1,q_1)=(0,0)$, or
                                            \item $s$ is odd and  $(p_1,q_1)=(1,0)$, or
                                            \item $s$ is even and  $(p_1,q_1)=(0,1)$.
                                          \end{itemize}
  \item  If $\beta$ is reducible, then $r=1$, $\alpha$ is rigid and $(p_1,q_1)\in \left\{(k,0), (0,k) \ : \ k>0\right\}$.
  \item If $\beta$ is pseudo-Anosov, $\alpha$ is rigid distinct from above.
\end{itemize}
\end{prop}

The conjugacy problem in $B_n$ can be solved by building the finite set of braids with minimal canonical length in the conjugacy class, called super summit set \cite{EM}. Notice that the non-periodic Murasugi representatives are rigid. Hence all Murasugi representatives have minimal canonical length in their conjugacy class, so they belong to their super summit set. Conversely, any braid $\beta$ in a super summit set satisfies that either $\beta$ or $\tau(\beta)$ is a Murasugi representative.

\begin{lema}{\rm \cite{Mat-Bert}}\label{Matt-Bert}
Let $\beta$ be a braid in $B_3$ with canonical length at least $2$. Then the super summit set of $\beta$, $SSS(\beta)$, is the set of
rigid conjugates of $\beta$. Actually, $SSS(\beta)$ consists of either two closed orbits under cycling, conjugate to each other by $\Delta$, or one closed orbit under cycling, self conjugated by $\Delta$.
\end{lema}

\subsection{Train tracks of braids in $B_3$}

\begin{figure}[h!]
\centerline{\includegraphics*[scale=0.65]{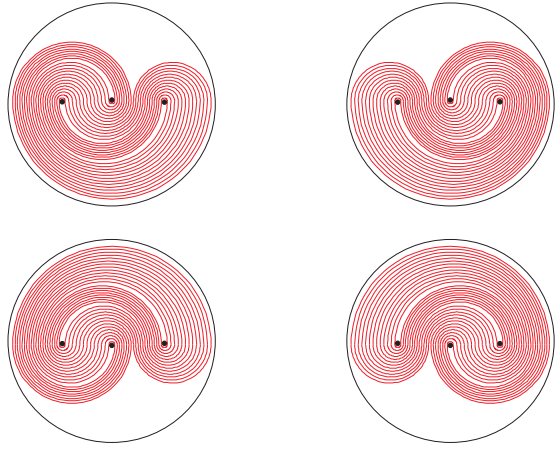}}
\caption{Foliations in $B_3$.}
\label{foliaions}
\end{figure}

Firstly, we want to point out the difference between train tracks as they have been defined in the section above, from those graphs that carry foliations. A train-track graph for a braid $\beta$ and the corresponding combinatorial map provide a matrix $M$. This matrix $M$ determines the Nielsen-Thurston type of $\beta$, and if it is pseudo-Anosov it determines the unstable foliation. However the reciprocal does not hold: from the foliation one cannot obtain a train track for $\beta$. Actually every admissible foliation in $\mathbb{D}_3$ (see Figure \ref{foliaions}) is \textit{carried} by one of the two graphs $\Gamma_1$ or $\Gamma_2$ in Figure \ref{2tt} \footnote{In literature, the graphs $G_1$, $G_2$ appear more often than $\Gamma_1$. The two foliations on the first row in Figure \ref{foliaions} are carried by $G_1$ and $G_2$ respectively.  We will use $\Gamma_1$ instead of $G_i$, despite the fact that it is not so intuitive to see that $\Gamma_1$ also carries both. Notice that if $x>y$ we could split $\Gamma_1$ to get $G_1$, changing the labels $\left\{ \begin{array}{l}
                                                                                  v_1=y\\
                                                                                  u_1=x-y>0
                                                                                \end{array}\right.$. And if $x<y$ we can split $\Gamma_1$ to get $G_2$. Similarly, $\Gamma_2$ can be split to get graphs $G_3$ and $G_4$, mirror images of $G_1$ and $G_2$.
\[\quad\]
\centerline{\includegraphics*[scale=0.9]{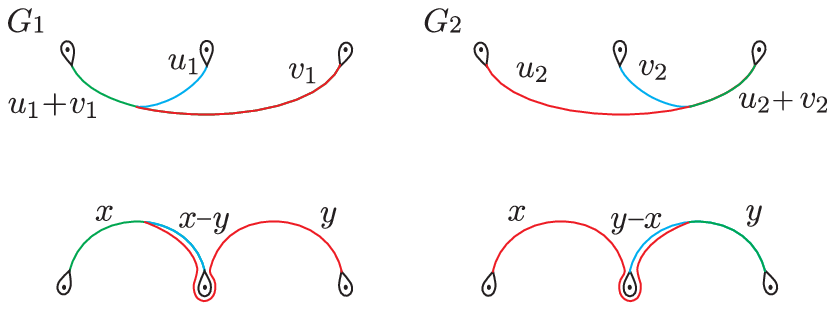}}}
with certain labels \cite{F-M}. However not every pseudo-Anosov braid admits $\Gamma_1$ or $\Gamma_2$ as train track, see the example below. Later it will be shown that if the braid is rigid, then yes, it admits either $\Gamma_1$ or $\Gamma_2$ as train track.

\begin{eje}
Neither $\Gamma_1$ nor $\Gamma_2$ are train tracks for the braid $\beta=(\sigma_1\sigma_1)^{-1}\sigma_2\sigma_1^{-1}(\sigma_1\sigma_1)$. The image of $\Gamma_i$ under the action of $\beta$ can not be embedded into a tubular neighborhood of $\Gamma_i$ respecting tangencies (see Figure \ref{ejemplo}).
\begin{figure}[h!]
\centerline{\includegraphics*[scale=1.1]{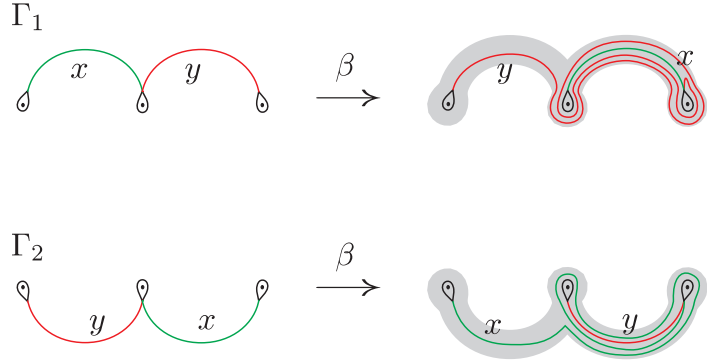}}
\caption{$\Gamma_1$, $\Gamma_2$ are not train track graph for $\beta=(\sigma_1\sigma_1)^{-1}\sigma_2\sigma_1^{-1}(\sigma_1\sigma_1)$.}
\label{ejemplo}
\end{figure}

We could also check the combinatorial maps $\widetilde{\beta}$ for $\Gamma_1$ and $\Gamma_2$. Let's label the edges around the punctures from left to right $e_1, e_2, e_3$. We will consider $e_i$ ($i=1,2,3$) oriented anticlockwise, and $x,y$ from right to left.
\[\begin{array}{rcl}
  \Gamma_1 & \stackrel{\widetilde{\beta}}{\rightarrow} & \Gamma_1 \\
  x & \rightarrow & y \\
  y & \rightarrow &(-x) \, e_2 \, (-y) \, e_3 \, y \, e_2 \, (-y) \, e_3
\end{array} \qquad \qquad
\begin{array}{rcl}
  \Gamma_2 & \stackrel{\widetilde{\beta}}{\rightarrow} & \Gamma_2 \\
  x & \rightarrow & y \\
  y & \rightarrow &(-y) \, (-x) \, e_3 \, x \, e_2 \,(-x)
\end{array}\]
The reader can check that in either case, $\widetilde{\beta}$ backtracks in the second iteration.
\end{eje}

The following theorem states that for a pseudo-Anosov rigid braid it is straightforward to get a train track and its transition matrix $M$. Actually, $M$ is the product of some simple matrices associated to the simple elements of the normal form. This is not the case in general: the product of transition matrices of two braids does not give any information about the product of the braids.

\begin{teo}{\rm \cite{Han}}
Let $\beta=(2s;1;p_1,q_1, \ldots, p_r, q_r)$ or $\beta=(2s+1;2;p_1,q_1, \ldots, p_r, q_r)$ be a rigid pseudo-Anosov braid in $B_3$. Then, the graph $\Gamma_1$ (see Figure \ref{2tt}) is a train-track graph for $\beta$.

Its transition matrix  is $M(\beta,\Gamma_1)=L^{p_1}U^{q_1}\cdots L^{p_r}U^{q_r}$, where
 \[L=\left(
     \begin{array}{cc}
       1 & 0 \\
       1 & 1 \\
     \end{array}
   \right) \qquad U=\left(
                    \begin{array}{cc}
                      1 & 1 \\
                      0 & 1 \\
                    \end{array}
                  \right).
 \]
\end{teo}

In fact this theorem can be easily extended to reducible braids. The graph $\Gamma_1$ is also a train tack for the braids $(2s;1;k,0)=\Delta^{2s}\sigma_1^{k}$, where $k>0$. Its associated matrix is $L^k=\left(
                                                                                   \begin{array}{cc}
                                                                                     1 & 0 \\
                                                                                     k & 1 \\
                                                                                   \end{array}
                                                                                 \right)$, which is obviously reducible.

If a braid $\beta$ has $\Gamma_1$ as train track and $M$ is the associated matrix, then $\tau(\beta)$ has $\Gamma_2=\Delta(\Gamma_1)$ as a train track. We order the edges by taking the one labeled by $x$ first, so $\beta$ and $\tau(\beta)$ have the same associated matrix $M$. Therefore the above theorem can be written for any rigid braid:

\begin{coro}
Every rigid braid in $B_3$ admits either $\Gamma_1$ or $\Gamma_2$ of Figure \ref{2tt} as a train track graph.
\end{coro}

  Let $\varepsilon$ be the automorphisms of  $B_n$ that maps $\sigma_i$ to $\sigma_i^{-1}$. For any braid $\beta$, $\tau\circ \varepsilon(\beta)$ is its vertical mirror image. This implies that $\tau\circ \varepsilon(\beta)$ has the same train track graph but it switches labels, hence both rows and columns of the transition matrix exchange places, $M(\tau\circ\varepsilon(\beta))=M(\varepsilon(\beta))=\left(                                                                                                                             \begin{array}{cc}
m_{2,2} & m_{2,1} \\
m_{1,2}& m_{1,1} \\                                                                                                                             \end{array}\right)$.
So, $\beta$ and $\varepsilon(\beta)$ have the same dilatation factor, and if $\overrightarrow{v}=(v_1, v_2)$ is the $M(\beta)$ eigenvector associated to $\lambda$, then $\overrightarrow{w}=(v_2,v_1)$ is the $M(\varepsilon(\beta))$ eigenvector associated to $\lambda$. This verifies an obvious fact: $\beta$ and $\varepsilon(\beta)$ have analogous dynamics.

Notice that the incidence matrix being a product of $L$'s and $U$'s, belongs to $SL_2(\mathbb{N})$. Furthermore, $L$ and $U$ are very special elements of this monoid. The following result, which is well known, is implicitly shown in the subsequent discussion.

\begin{teo}
The monoid $SL_2(\mathbb{N})$ is freely generated by $L$, $U$.
\end{teo}

This theorem means that given any $2\times 2$ matrix $M$ with non negative entries and determinant $1$, it admits a unique product decomposition in terms of $L$ and $U$. We will call this the $LU$-decomposition of $M$. Because $det(M)=1$ and it has nonzero entries, we can say that $M$ has a biggest row $R^{(M)}$ and a smallest row $r^{(M)}$, meaning $R^{(M)}_1\geq r^{(M)}_1$ and $R^{(M)}_2\ge r^{(M)}_2$ (where at least one of the inequalities is strict). If $r^{(M)}=\left(m_{1,1}, m_{1,2}\right)$, then the first factor of the $LU$-decomposition of $M$ is $L$:
\[M=\left(
      \begin{array}{c}
        r^{(M)} \\
        R^{(M)} \\
      \end{array}
    \right)\ \Rightarrow \ M=LM', \ \mbox{where } M'= \left(
                             \begin{array}{c}
                               r^{(M)} \\
                               R^{(M)}-r^{(M)} \\
                             \end{array}
                           \right)\in SL_2(\mathbb{N}).
\]
And obviously, if the first row is $R^{(M)}$, then the first factor is $U$. This gives us an algorithm to compute the $LU$-decomposition of $M$. At each step $\|M\|_{\infty}>\|M'\|_{\infty}$, therefore the algorithm ends. We could have defined the analogous algorithm defining the columns $C^{(M)}$ and $c^{(M)}$, obtaining the decomposition from right to left.

Given $M\in SL_2(\mathbb{N})$, we define the length of $M$, and we write $\ell(M)$, as the length on the associated $LU$-decomposition. We will denote $M(i)$ the product of the first $i$ factors of $M$, for $i=1,\ldots,\ell(M)$. With this notation $M(\ell(M))=M$, the first factor of $M$ is $M(1)$ and the $i^{\mbox{th}}$ factor is $M(i-1)^{-1}M(i)$.

\subsection{Rigid $3$-braids and $SL_2(\mathbb{N})$}

 As a consequence of the results in the section above, we can prove Theorem \ref{corointro}.\\

\begin{demo}[of theorem\ref{corointro}]
The map from rigid braids to matrices in $SL_2(\mathbb{N})$ is explicit in both senses. Given the normal form of a rigid braid $\beta=\Delta^{s}s_1\cdots s_k$, we can construct from the associated tuple $(s;i;p_1,q_1, \ldots, p_r, q_r)$ the $LU$-word $L^{p_1}U^{q_1}\cdots L^{p_r}U^{q_r}$. Reciprocally, given a $2\times2$ square matrix $M$, the $LU$-decomposition provides the values $p_j, q_j$. The value $s$ must be chosen even or odd, so that the associated braid is rigid. The $2\times2$ identity matrix is associated to powers of $\Delta$, so $1$ and $\Delta$ (mod $\Delta^2$) are the two preimages of the identity matrix. In the other cases, because there are always two options for $i$, the correspondence restricted to $B_n/\langle \Delta^2 \rangle$ is 2 to 1.
\end{demo}

Due to Lemma \ref{Matt-Bert}, it follows that two rigid conjugated braids must have conjugate matrices in $SL_2(\mathbb{N})$, that is cyclic permutations of $L$, $U$ factors.

Let $J=\left(
                                                      \begin{array}{cc}
                                                        0 & 1 \\
                                                        1 & 0 \\
                                                      \end{array}
                                                    \right)$. Because $J=J^{-1}$ and $L=JUJ$, we have that $M(\varepsilon(\beta))=JM(\beta)J$. The $LU$-decomposition of a matrix $M$ can be obtained from that of $JMJ$ by exchanging $L$'s and $U$'s. However, these two matrices are not conjugated in $SL(\mathbb{N})$ in general.

Due to the above results, we can explicitly give the set of dilatation factors in $B_3$:

\begin{coro}
The set of dilatation factors for pseudo-Anosov braids in three strands is:
\[ \left\{ \lambda=\frac{T+ \sqrt{T^2-4}}{2}, \ T\in \mathbb{N}, T\geq 3 \right\}.\]
\end{coro}

\begin{demo}
If the incident matrix associated to a braid is $M\in SL_2(\mathbb{N})$, the dilatation factor is the biggest root of the polynomial $x^2-Tx+1$, where $T=trace(M)\in \mathbb{N}$. So the dilatation factor $\lambda$ of a pseudo-Anosov braid only depends on $T=trace(M)$. As $det(M)=1$ and all entries of $M$ are non-negatives integers, we have $T>1$. Only $L^k, U^k$ have trace equal to $2$, but these matrices are reducible. Therefore, for pseudo-Anosov braids $T>2$, and every value of $T>2$ can be obtained at least for the matrix $LU^{T-2}=\left(
                                                                \begin{array}{cc}
                                                                  1 & T-2 \\
                                                                  1 & T-1 \\
                                                                \end{array}
                                                              \right)$, which corresponds, among others, to the braid $\sigma_1\sigma_2^{-(T-2)}$.

\end{demo}

\subsection{Canonical length and dilatation factor}

In this subsection we will finally relate the canonical length of a rigid braid with its dilatation factor. We remark that a rigid braid has minimal canonical length in its conjugacy class, hence these results will not only provide a lower bound for the dilatation factor  of a rigid pseudo-Anosov $3$-braid with fixed canonical length, but also an upper bound for the dilatation factor  of any pseudo-Anosov $3$-braid with fixed canonical length.

As the canonical length of a rigid braid is precisely the length of its $LU$-decomposition, and the dilatation factor depends only on the trace of the associated matrix, we just need to relate the trace of a matrix in $SL_2(\mathbb{N})$ with its $LU$-decomposition length $\ell(M)$.

 \begin{nota}
 These two conjugacy class invariants, dilatation factor and minimal length, do not characterize the conjugacy class of a braid. We give as example the rigid braids $\beta_1$ and $\beta_2$, which are not conjugate\footnote{Due to the homogeneity of the relations in $B_3$, if two braids are conjugate, the sum of the exponents of generators in their representing words must be the same. Hence, modulo $\Delta^2$ the sums must be the same modulo $6$.} modulo $\Delta^2$:
\[\begin{array}{lll}
    \beta_1= \Delta \sigma_1^3 \ \sigma_1\sigma_2 \  \sigma_2 \ \sigma_2\sigma_1 \ \sigma_1 \ \sigma_1\sigma_2 & \quad & M(\beta_1)=L^3 U L U L U = \left(
                                                                                                                       \begin{array}{cc}
                                                                                                                         5 & 8 \\
                                                                                                                         18 & 29 \\
                                                                                                                       \end{array}
                                                                                                                     \right)\\
    \beta_2=\left(\sigma_1^2\ \sigma_1\sigma_2 \ \sigma_2\sigma_1\right)^2 & \quad &  M(\beta_2)=L^2 U^2 L^2 U^2 = \left(
                                                                                                                       \begin{array}{cc}
                                                                                                                         5 & 12 \\
                                                                                                                         12 & 29 \\
                                                                                                                       \end{array}
                                                                                                                     \right)
  \end{array}\]
\end{nota}

Now, we give a lower bound for the trace of a matrix with fixed length.

\begin{prop}
For a given length $\ell$, the irreducible matrix in $SL_2(\mathbb{N})$ with lowest trace, up to conjugation in $SL_2(\mathbb{N})$ and conjugation by $J$, is:
\[LU^{\ell-1}=\left(
             \begin{array}{cc}
               1 & \ell-1 \\
               1 & \ell \\
             \end{array}
           \right)
\]
\end{prop}

\begin{demo}
Let $M$ be an irreducible matrix with minimal trace and $\ell(M)=\ell$. As $M$ is irreducible, $\ell>1$ and, up to conjugation in $SL_2(\mathbb{N})$ we can suppose that $M$ starts with $LU=\left(
                                                                                                                             \begin{array}{cc}
                                                                                                                               1 & 1 \\
                                                                                                                               1 & 2 \\
                                                                                                                             \end{array}
                                                                                                                           \right)$.
For any $2\times2$ matrix $N=\left(
                        \begin{array}{c|c}
                          C_1 & C_2 \\
                        \end{array}
                      \right)$, where $C_1,\ C_2$ are column vectors,
                      \[NL=\left(\begin{array}{c|c}
                                    C_1+C_2 & C_2 \\
                    \end{array} \right) \qquad NU=\left(
                                              \begin{array}{c|c}                                                                                               C_1 & C_1+C_2 \\
                                              \end{array}\right).\]
Therefore $trace(NL)=trace(N)+n_{1,2}$ and $trace(NU)=trace(N)+n_{2,1}$. Setting $N=LU$ which has positive entries, at each multiplication by $L$ or $U$ the trace increases at least one unit: hence $trace\left(M(i)\right)\geq i+1$ for $i=1,\ldots,\ell$. As $\ell(M)=\ell$, it follows that $trace(M)\geq \ell +1$.

If $\ell>2$, the third factor of $M$ can be either $U$ or $L$. In the first case, as $LU^2=\left(
                                                                                             \begin{array}{cc}
                                                                                               1 & 2 \\
                                                                                               1 & 3 \\
                                                                                             \end{array}
                                                                                           \right)$, every subsequent factor of $M$ must also be $U$, otherwise the trace would increase more than required. Hence $M=LU^{\ell-1}$.
In the second case, as $LUL=\left(
                              \begin{array}{cc}
                                2 & 1 \\
                                3 & 2 \\
                              \end{array}
                            \right)$, every subsequent factor must be $L$. Hence $M=LUL^{\ell-2}$, which is conjugate to $UL^{\ell-1}$ in $SL_2(\mathbb{N})$, which is conjugate by $J$ to $LU^{\ell-1}$.
\end{demo}

From the above, the following result is straightforward.
\begin{coro}\label{coro}
Given $M\in SL_2(\mathbb{N})$, $\ell(M) \leq trace(M)-1$.
\end{coro}

Now we give an upper bound for the trace.

\begin{prop}\label{upperb}
For a given length $\ell>2$, the matrix in $SL_2(\mathbb{N})$ with greatest trace, up to conjugation in $SL_2(\mathbb{N})$ and conjugation by $J$, is:
\[\begin{array}{lr}
  \left(LU\right)^{\frac{\ell}{2}}=\left(
                  \begin{array}{cc}
                    F_{\ell-1} & F_\ell \\
                    F_\ell & F_{\ell +1} \\
                  \end{array}
                \right)
   & \mbox{if } \ell \mbox{ is even,} \\
  \left(LU\right)^{\frac{\ell-1}{2}}L=\left(
                  \begin{array}{cc}
                    F_{\ell} & F_{\ell-1} \\
                    F_{\ell+1} & F_{\ell} \\
                  \end{array}
                \right)
   & \mbox{if } \ell \mbox{ is odd,}
\end{array}\]
where $F_0=0, \ F_1=1, \ F_{i+2}= F_{i+1}+ F_i$. That is, $F_\ell$ is the $\ell^{th}$ Fibonacci number.
\end{prop}

\begin{demo}
Let $A=LULU\cdots$ be the matrix of length $\ell$ defined in the statement, and let $B$ be another matrix with the same length. If $B$ is reducible, its trace must be $2$, which is lower than $trace(A)$. If $B$ is irreducible, up to cyclic permutation of the factors and conjugation by $J$, we can assume that both $LU$-decompositions of $A,\ B$ start and finish with the same matrices.

   Let us prove by induction that $C^{(A(i))}\geq C^{(B(i))}$, and $c^{(A(i))}\geq c^{(B(i))}$, $i=1,\ldots, \ell$, where $C^{(N)}$ (resp. $c^{(N)}$) is the biggest column of a matrix $N\in SL_2(\mathbb{N})$.
   It is easy to see that we construct $A(i+1)$ keeping the biggest column of $A(i)$ and replacing the smallest by the sum of both columns in $A(i)$:
  \[ \begin{array}{l}
        C^{(A(i+1))}=C^{(A(i))}+c^{(A(i))} \\
       c^{(A(i+1))}=C^{(A(i))}.
     \end{array}\]
  For $B(i+1)$, all we can say is
  \[  \begin{array}{l}
      C^{(B(i+1))}=C^{(B(i))}+c^{(B(i))} \\
      c^{(B(i+1))}\leq C^{(B(i))}.
  \end{array}\]
  We have $A(1)=B(1)=L$, so the claim holds for $i=1$. If $C^{(A(i))}\geq C^{(B(i))}$, and $c^{(A(i))}\geq c^{(B(i))}$ one has:
  \[\begin{array}{l}
      C^{(A(i+1))}=C^{(A(i))}+c^{(A(i))}\geq C^{(B(i))}+c^{(B(i))}=C^{(B(i+1))} \\
      c^{(A(i+1))}=C^{(A(i))}\geq C^{(B(i))} \geq c^{(B(i+1))}
    \end{array}\]
and the claim is shown.

Now, as the $\ell^{\mbox{th}}$ factors in the $LU$-decompositions of $A$ and $B$ are the same, both $C^{(A)}$ and $C^{(B)}$ are in the same relative position: either both of them are the first columns of the corresponding matrices or they are the second ones. This implies $tr(A)\geq tr(B)$, and the equality holds if and only if $A=B$.

We will prove by induction that the entries of the product $LULU\ldots$ are Fibonacci numbers as stated. If $\ell=1$, then $L=\left(
                                                                                                                                     \begin{array}{cc}
                                                                                                                                       F_1 &F_0 \\
                                                                                                                                        F_2& F_1 \\
                                                                                                                                     \end{array}
                                                                                                                                   \right)$. So if the claim holds for odd length $\ell$, then we have:
\[\begin{array}{l}
  LU\cdots LUL * U=\left(
                    \begin{array}{cc}
                        F_{\ell} & F_{\ell-1} \\
                        F_{\ell +1} & F_{\ell} \\
                    \end{array}\right) \left(
                                         \begin{array}{cc}
                                           1 & 1 \\
                                           0 & 1 \\
                                         \end{array}
                                       \right)=\left(
                                                         \begin{array}{cc}
                                                           F_{\ell} & F_{\ell+1} \\
                                                           F_{\ell+1} & F_{\ell+2} \\
                                                         \end{array}
                                                       \right)\\
    LU\cdots LU * L=\left(
                    \begin{array}{cc}
                      F_{\ell-1} & F_\ell \\
                    F_\ell & F_{\ell +1} \\
                    \end{array}
                  \right)\left(
                           \begin{array}{cc}
                             1 & 0 \\
                             1 & 1 \\
                           \end{array}
                         \right)=\left(
                                           \begin{array}{cc}
                                             F_{\ell +1} & F_\ell \\
                                             F_{\ell +2} & F_{\ell +1} \\
                                           \end{array}
                                         \right)
                                          \end{array}\]

\end{demo}

Now we can prove our main result.\\

\begin{demo}[of Theorem \ref{teo1}]
Let $\beta\in B_3$ be a pseudo-Anosov rigid braid. Its Garside length $\ell$ must be at least $2$. Then $\beta$ has dilatation factor $\lambda=\frac{1}{2}\left( T+ \sqrt{T^2-4}\right)$, for an integer $T\geq 3$. The lower bound for $\lambda$ is direct from the fact that $\ell+1 \leq T$ by Corollary \ref{coro}.

Let $\lambda(M)$ the greatest eigenvalue of a matrix $M\in SL_2(\mathbb{N})$. As a consequence of Proposition \ref{upperb}, for the upper bound when $\ell$ is even, we have that
\[\lambda\left(\left(LU\right)^{\frac{\ell}{2}}\right)=\lambda\left(LU\right)^{\frac{\ell}{2}}=\left(1+\phi\right)^{\frac{\ell}{2}}=\phi^{\ell}.\]

 If $\ell$ is odd, we get that the greatest dilatation factor of a braid with odd canonical length is $\lambda=\frac{1}{2}\left( 2F_{\ell}+ \sqrt{\left(2F_{\ell}\right)^2-4}\right)=F_{\ell}+\sqrt{F_{\ell}^2-1}$. Therefore we have that
\[2F_{\ell}-1<\lambda < 2F_{\ell}.\]
  As $F_n=\frac{\phi^n-(-\phi)^{-n}}{\sqrt{5}}$, for the case $\ell$ odd we have that
  \[\frac{2}{\sqrt{5}}\left( \phi^\ell+\phi^{-\ell}\right) -1 < \lambda < \frac{2}{\sqrt{5}}\left( \phi^\ell+\phi^{-\ell}\right).\]
 \end{demo}

\end{document}